# Performance evaluation through DEA benchmarking adjusted to goals

*José L. Ruiz and Inmaculada Sirvent*

*Centro de Investigación Operativa. Universidad Miguel Hernández. Avd. de la Universidad, s/n 03202-Elche (Alicante), SPAIN*

*July, 2018*

## Abstract

Data Envelopment Analysis (DEA) is extended to the evaluation of performance of organizations within the framework of the implementation of plans for improvements that set management goals. Managers usually set goals without having any evidence that they will be achievable at the moment of conducting performance evaluation or, on the contrary, they may set little too unambitious goals. Using DEA for the benchmarking ensures an evaluation in terms of targets that both are attainable and represent best practices. In addition, the approach we propose adjusts the DEA benchmarking to the goals in order to consider the policy of improvements that was pursued with the setting of such goals. From the methodological point of view, the models that minimize the distance to the DEA strong efficient frontier are extended to incorporate goal information. Specifically, the models developed seek DEA targets that are as close as possible to both actual performances and management goals. To illustrate, we examine an example that is concerned with the evaluation of performance of public Spanish universities.

## 1. Introduction

In management, organizations often set goals in planning improvements. Supervisors and managers set the goals in a given period of time and an ex post evaluation of performance relative to such goals is carried out as part of a process of monitoring and control. Goals are typically established as the result of interactions between stakeholders, taking into account objectives, policies, the knowledge of prior period performance, etc. An example of this situation can be found in the public Spanish universities, wherein regional governments and university managers formalize an agreement (a programme-contract) in which goal levels are set for a number of index indicators related to the different areas of performance (teaching, research, knowledge transfer,…). In fact, part of the financing of public universities is linked to the achievement of goals, so that a payment of incentives is made in the next budgetary year on the basis of an evaluation of performance relative to the goals.

The focus in the present paper is on the evaluation of performance in organizations within the framework of the implementation of plans for improvements that set management goals. In respect of the setting of goals, it should be highlighted, on one hand, that stakeholders may set goals without having any evidence that they will be achievable. And on the other, that managers whose own performance is to be evaluated often participate in the setting of goals, so that they have the opportunity to influence by setting little too unambitious goals, while improvements should be directed towards best practices. Taking into account these considerations, we develop an approach that make it possible an evaluation of performance based on targets that are attainable as well as represent best practices.

In order to do so, we propose the use of Data Envelopment Analysis (DEA) (Charnes et al., 1978). Cooper (2005) emphasizes DEA as a tool *directed to evaluating past performances as part of the control function of management*. Specifically, DEA evaluates performance of decision making units (DMUs) from a perspective of benchmarking, through setting targets on the efficient frontier of a production possibility set (PPS). We highlight the fact that, as stated in Cook et al. (2014), *In the circumstance of benchmarking, the efficient DMUs, as defined by DEA, may not necessarily form a "production frontier", but rather lead to a "best-practice frontier"*. In the context of the implementation of plans for improvements, we develop a DEA-based benchmarking approach that allows us to set targets that represent best practices, while at the same time considering what is technically achievable at the moment of conducting performance evaluation. Note that, while goals have been previously established and exogenously, and so they may not be either achievable or efficient, DEA targets both are attainable and represent best practices, insofar as they are set from actual performances (of the evaluation period) through benchmark selection among points on to the DEA strong efficient frontier of the PPS. In addition to setting targets achievable and efficient, the proposed approach adjusts the DEA benchmarking to the goals. Relating DEA targets and goals seeks to consider the policy of improvements that was pursued through the setting of goals (note that DMUs might have oriented their activities towards the achievement of the goals). Concretely, the models formulated set DEA targets that are as close as possible to both actual performances and goals. Thus, on one hand, targets are related to actual performances, in order to allow for the individual circumstances of the DMUs under evaluation and, on the other, the benchmarking is adjusted in order to evaluate performance not only considering best practices (technically achievable) but also what is desirable, as represented by the goals.

In Stewart (2010) target setting is addressed within a general framework that incorporates long term goals to the DEA models. Specifically, the author finds realistically achievable targets on the efficient frontier through a search that is made taking into account the goals (like us here) by using the Wierzbicki reference point approach. Nevertheless, in Stewart (2010) the aim is to develop an ex ante planning tool while we are concerned with ex post evaluations within a monitoring process of improvements. Following Stewart (2010), goal directed benchmarking models are proposed in Azadi et al. (2013) and Azadi et al. (2014), for supplier selection, and in Khoveyni and Eslami (2016), for the case of having some inputs that are imposed on the DMUs. Cook et al. (2018) have recently proposed a DEA-based benchmarking approach in the context of pay-for-performance incentive plans. Target setting there is also carried out on the DEA efficient frontier, by minimizing the differences between the payments relative to the goals and those relative to the targets. Thanassoulis and Dyson (1992) develop a model that estimates targets on the basis of "ideal" target input-output levels the DMUs may be able to specify. Since these levels may be neither feasible nor efficient, a two-stage process is followed in order to determine a set of efficient input-output levels which are compatible with the "ideal" targets. Thanassoulis and Dyson (1992) actually develop several benchmarking models that incorporate a priori preference specifications, albeit such preferences are not exactly expressed as goals in terms of input-output levels but through the specification of some user-supplied weights that represent the relative desirability of improving the levels of the different inputs and outputs. Other authors have extended these models: for example, Zhu (1996) relaxes the restrictions on the factors that represent the proportions by which inputs and outputs are changed while Lins et al. (2004) develop models with a posteriori preference specifications with the purpose of alleviating the problems related to the specification of the weights. In Athanassopoulos (1995) extensions of DEA to the global organizational level that embark from Thanassoulis and Dyson (1992) are proposed, through models that allow to accommodating global organizational targets (see also Athanassopoulos (1998)). For other models allowing for a priori preferences of decision makers see Lozano and Villa (2005), Hinojosa and Mármol (2011) and Fang (2015). Besides, the DEA models that incorporate judgmental preferences through weight restrictions (Thompson et al. (1986), Charnes et al. (1990)) can also be considered in this class of models (see Ramón et al. (2016) for a discussion on the use of DEA models with weight restrictions for benchmarking and target setting). An alternative approach is to incorporate preferences in an interactive fashion, without requiring prior judgments. The models developed under this approach usually combine DEA and multi objective linear programming techniques. Lins et al. (2004), mentioned above, also includes an interactive multiple objective target optimization. See also Golany (1988), Post and Spronk (1999), Halme et al. (1999), Joro et al. (2003) and Yang et al. (2009).

Methodologically, the DEA models that minimize the distance to the efficient frontier are extended here to incorporate information on goals. The models that set the closest targets have made an important contribution to DEA as tool for the benchmarking. Closest targets minimize the gap between actual performances and best practices, thus showing the DMUs the way for improvement with as little effort as possible. In Aparicio et al. (2007), the problem of minimizing the distance to the efficient frontier was theoretically solved through a mixed integer linear problem (MILP) that can be used to set the closest targets. The methodology introduced in that paper has been extended in the last decade to deal with different issues related to target setting and benchmarking: in Ramón et al. (2016) to deal with models including weight restrictions, in Ruiz and Sirvent (2016) for developing a common benchmarking framework, in Aparicio et al. (2017) it is extended to oriented models, in Cook et al. (2017) for the benchmarking of DMUs classified in groups, in Ramón et al. (2018), which propose a sequential approach for the benchmarking and in Cook et al. (2018), above mentioned, for target setting in pay-for performance incentive plans (see Portela et al. (2003) and Tone (2010) for other papers on closest targets in DEA)[1]. Within this methodological framework, we now propose new benchmarking models that incorporate management goals. Specifically, we develop both non-oriented and oriented models, while taking into account the possibility of having non-controllable factors that need to be considered for performance evaluation. These models set targets on the DEA strong efficient frontier by considering two objectives in their formulations: minimizing the distance between actual inputs and/or outputs and targets and minimizing the distance between targets and goals. Thus, as has been said before, DEA targets in this approach allow managers to evaluate the degree of improvements taking into consideration as much as possible both actual performances and goals. To illustrate, we examine an example about teaching performance of public Spanish universities in the academic year 2014-15, by assuming that some goals for several outputs regarding the progress, retention and graduation of students were set in 2013-14, and considering at the same time their use of resources, in terms of academic staff and expenditures, as well as the size of universities as a non-controllable factor.

The paper unfolds as follows: In section 2 we develop a general approach for setting targets through DEA benchmarking adjusted to goals. Section 3 addresses the situation where goals have been established only for inputs or for outputs, by using DEA oriented models. In section 4 we adapt the models developed for the case of having inputs/outputs whose levels are exogenously fixed. Section 5 includes an empirical illustration. Last section concludes.

## 2. Target setting through DEA benchmarking adjusted to goals

Throughout this section we consider that we have *n* DMUs whose performances in a given period of time are to be evaluated in terms of *m* inputs and *s* outputs. These are denoted by $\left(X_j, Y_j\right)$, $j=1,...,n$, where $X_j = \left(x_{1j},...,x_{mj}\right)' > 0_m$, $j=1,...,n$, and $Y_j = \left(y_{1j},...,y_{sj}\right)' > 0_s$, $j=1,...,n$.[2]

---

[1] We note that the idea of closeness to the efficient frontier has also been investigated for purposes of developing efficiency measures that satisfy some desirable properties (see Ando et al. (2017), Aparicio and Pastor (2014), Fukuyama et al. (2014a, 2014b) and Fukuyama et al. (2016)).

We suppose that, at a prior stage to the evaluation of performance, some goals were established within the framework of the implementation of a plan for improvements. These goals represent levels of the inputs and outputs the DMUs aspire to achieve in the subsequent performance period t. For every DMU$_j$, these goals are denoted by $G_j^I = \left(g_{1j}^I,....,g_{mj}^I\right)'$, for the inputs, and $G_j^O = \left(g_{1j}^O,....,g_{sj}^O\right)'$, for the outputs. Performance is evaluated through a benchmarking analysis, by assuming a variable returns to scale (VRS) technology (Banker et al., 1984). Thus, the production possibility set (PPS), $T = \{(X,Y) / X \text{ can produce } Y\}$, can be characterized as $T = \left\{(X,Y) \middle/ X \geq \sum_{j=1}^{n} \lambda_j X_j, \ 0 \leq Y \leq \sum_{j=1}^{n} \lambda_j Y_j \right.$

$\left. \sum_{j=1}^{n} \lambda_j = 1, \lambda_j \geq 0 \right\}$ (the approach proposed can be similarly developed in terms of a constant returns to scale (CRS) technology).

In this situation, the evaluation of performances can be seen as part of a monitoring process of improvements. In order to carry out such evaluation, we propose to use DEA, which allows us to set targets that establish the potential for improvement in inputs and outputs for the DMUs. In general, using a DEA approach makes it possible to compare actual performances of DMUs against benchmarks on the efficient frontier, thus ensuring an evaluation in terms of targets that both are attainable and represent best practices of their peers. In particular, a DEA analysis based on closest targets shows the DMUs in their best possible light. However, when information on goals is available, it can be argued that the benchmarking should be adjusted in order to relate DEA targets and management goals, aside from relating DEA targets and actual performances. That is, benchmarking models should consider the strategic planning established by stakeholders in addition to considering the individual circumstances of the DMUs. For this reason, we need to formulate a DEA model that provides an evaluation of performance based on targets that satisfy the following two objectives (as much as possible):

1. DEA targets should be as close as possible to actual performances.
2. DEA targets should be as close as possible to management goals.

Formally, we will find the DEA targets for a given DMU$_0$ that are wanted by trying to minimize the two deviations below (expressed in terms of a weighted L$_1$-norm)

---

[2] Positive inputs and outputs are required because we use a weighted L$_1$-norm in which deviations between targets and both actual performances and goals are calculated relative to actual inputs/outputs of the unit under evaluation.

$$\left\|(X_0,Y_0)-(\hat{X}_0,\hat{Y}_0)\right\|_1^\omega = \sum_{i=1}^{m}\frac{|x_{i0}-\hat{x}_{i0}|}{x_{i0}} + \sum_{r=1}^{s}\frac{|y_{r0}-\hat{y}_{r0}|}{y_{r0}} \qquad (1.1)$$

$$\left\|(G_0^I,G_0^O)-(\hat{X}_0,\hat{Y}_0)\right\|_1^\omega = \sum_{i=1}^{m}\frac{|g_{i0}-\hat{x}_{i0}|}{x_{i0}} + \sum_{r=1}^{s}\frac{|g_{r0}-\hat{y}_{r0}|}{y_{r0}} \qquad (1.2)$$

(1)

over the points $(\hat{X}_0,\hat{Y}_0)$ on the strong efficient frontier of the PPS, which is denoted by $\partial(T)$.

Ideally, one would like to minimize simultaneously both (1.1) and (1.2). However, it is obvious that the targets that minimize (1.1) will not be necessarily those that minimize (1.2), and the other way around too. As a compromise between these two objectives, we propose to minimize a convex combination of them

$$\alpha\left\|(X_0,Y_0)-(\hat{X}_0,\hat{Y}_0)\right\|_1^\omega + (1-\alpha)\left\|(G_0^I,G_0^O)-(\hat{X}_0,\hat{Y}_0)\right\|_1^\omega, \qquad (2)$$

where $0\leq\alpha\leq 1$, over $\partial(T)$. This will provide us with non-dominated solutions[3]. Through the specification of $\alpha$, we may adjust the importance that is attached to each of the two objectives. The specification $\alpha=1$ leads to the closest targets to DMU$_0$ on $\partial(T)$. Obviously, those targets would be established ignoring completely the goals. As $\alpha$ decreases, targets are more in line with the goals. At the extreme, i.e., when $\alpha=0$, we find the closest targets to the goals on $\partial(T)$. In that case, targets would be established ignoring completely the individual circumstances of DMU$_0$, so they may not necessarily be acceptable for the managers; this will depend on the goal levels that were set. In practice, it may be useful to provide a series of benchmarks that are generated by specifying values for $\alpha$ in a grid within [0,1] (this approach has already been used in Stewart (2010), wherein the specification of $\alpha$ determines different reference points). Managers may appreciate having available a spectrum of benchmarks which provides an overview of the environment, so that an evaluation of improvements can be made taking into account where they stand in the industry or sector. This all will be illustrated in section 5.

Bearing in mind the above, a model that provides the targets that are sought for DMU$_0$ (for a given $\alpha$) can be formulated as follows

---

[3] It is known that the weighted sum method to solve a multi-objective optimization problem will produce Pareto optimal solutions when all the weights are strictly positive. However, when some of the weights is zero, there is a potential to get only weakly efficient solutions. For this reason, to ensure non-dominated solutions in the cases α=0 and 1, a second stage is needed, in which the objective that is not considered in (2) is minimized subject to (2) takes the optimal value obtained in the first stage.

$$\text{Min} \quad \alpha \left\| (X_0, Y_0) - (\hat{X}_0, \hat{Y}_0) \right\|_1^\omega + (1-\alpha) \left\| (G_0^I, G_0^O) - (\hat{X}_0, \hat{Y}_0) \right\|_1^\omega$$
$$\text{s.t.:} \tag{3}$$
$$(\hat{X}_0, \hat{Y}_0) \in \partial(T)$$

The following proposition provides a characterization of $\partial(T)$.

**Proposition 1 (Ruiz et al., 2015).**

$(\hat{X}_0, \hat{Y}_0) \in \partial(T) \Leftrightarrow \hat{X}_0 = \sum_{j \in E} \lambda_j X_j, \quad \hat{Y}_0 = \sum_{j \in E} \lambda_j Y_j, \quad \sum_{j \in E} \lambda_j = 1,$ for some $\lambda_j \geq 0, \quad j \in E,$ satisfying $\lambda_j d_j = 0, \quad j \in E,$ where $-v'X_j + u'Y_j + u_0 + d_j = 0, \quad j \in E, \quad v \geq 1_m, \quad u \geq 1_s, \quad u_0 \in \mathbb{R},$ E being the set of extreme efficient units of T. [4,5]

Taking into account this proposition, we may derive the following operational formulation of model (3)

---

[4] See Charnes et al. (1991) for a classification of DMUs. In particular, identifying the DMUs in E can be done by simply checking if the optimal value of the following problem $\text{Min} \left\{ \lambda_0 \Big/ \sum_{j=1}^n \lambda_j x_{ij} \leq x_{i0}, i = 1, ..., m, \sum_{j=1}^n \lambda_j y_{rj} \geq y_{r0}, r = 1, ..., s, \sum_{j=1}^n \lambda_j = 1, \lambda_j \geq 0, j = 1, ..., n \right\}$ is 1.

[5] In Ruiz et al. (2015), the proposition is equivalently enunciated by using the classical big M and binary variables.

$$\text{Min} \quad \alpha\left(\sum_{i=1}^{m}\frac{|t_{i0}|}{x_{i0}}+\sum_{r=1}^{s}\frac{|s_{r0}|}{y_{r0}}\right)+(1-\alpha)\left(\sum_{i=1}^{m}\frac{|t_{i0}^{g}|}{x_{i0}}+\sum_{r=1}^{s}\frac{|s_{r0}^{g}|}{y_{r0}}\right)$$

s.t.:

$$\sum_{j\in E}\lambda_{j}x_{ij}=x_{i0}-t_{i0} \qquad i=1,...,m$$

$$\sum_{j\in E}\lambda_{j}y_{rj}=y_{r0}+s_{r0} \qquad r=1,...,s$$

$$\sum_{j\in E}\lambda_{j}x_{ij}=g_{i0}^{I}-t_{i0}^{g} \qquad i=1,...,m$$

$$\sum_{j\in E}\lambda_{j}y_{rj}=g_{r0}^{O}+s_{r0}^{g} \qquad r=1,...,s$$

$$\sum_{j\in E}\lambda_{j}=1$$

$$-\sum_{i=1}^{m}v_{i}x_{ij}+\sum_{r=1}^{s}u_{r}y_{rj}+u_{0}+d_{j}=0 \qquad j\in E$$

$$v_{i}\geq 1 \qquad i=1,...,m$$

$$u_{r}\geq 1 \qquad r=1,...,s$$

$$\lambda_{j}d_{j}=0 \qquad j\in E$$

$$d_{j},\lambda_{j}\geq 0 \qquad j\in E \qquad (4)$$

$$t_{i0},t_{i0}^{g} \text{ free} \qquad i=1,...,m$$

$$u_{0},s_{r0},s_{r0}^{g} \text{ free} \qquad r=1,...,s$$

The key in model (4) is in the constraints $\lambda_{j}d_{j}=0$, $j\in E$. By virtue of these, model (4) ensures targets on the Pareto efficient frontier of T, $\partial(T)$. Note that, if $\lambda_{j}>0$ then $d_{j}=0$. This means that if $DMU_{j}$ in E participates actively as a referent in the benchmarking of $DMU_{0}$ then it necessarily belongs to $-v'X+u'Y+u_{0}=0$. That is, the extreme efficient DMUs that participate in the benchmarking of the $DMU_{0}$ are all on a same facet of $\partial(T)$, because they all belong to the same supporting hyperplane of T, $-v'X+u'Y+u_{0}=0$, whose coefficients are non-zero. Therefore, solving (4) allows us to identify a reference set of efficient DMUs, $RS_{0}=\{DMU_{j}/\lambda_{j}^{*}>0\}$, for the benchmarking of $DMU_{0}$, so that targets result from its projection onto the facet of $\partial(T)$ spanned by the DMUs in $RS_{0}$. Specifically, these targets are actually the coordinates of the projection of $DMU_{0}$ on $\partial(T)$ that results of minimizing (2). In short, Proposition 1 allows us to handle directly projection points of $DMU_{0}$ on $\partial(T)$ as feasible solutions of the benchmarking models, so that targets result from minimizing deviations between potential benchmarks on $\partial(T)$ and both actual performances and goals. Thus, we address explicitly the two objectives that have been established (in Stewart (2010), projections of DMU$_0$ on $\partial(T)$ are found by minimizing deviational variables with respect to a reference point, which is a combination of actual performances and goals, following the Wierzbicki approach).

For every $\alpha$, the targets $\left(\hat{X}_0^*, \hat{Y}_0^*\right)$ provided by (4) can be expressed in terms of its optimal solutions as

$$\begin{aligned}\hat{X}_0^* &= \sum_{j \in E} \lambda_j^* X_j \left(= X_0 - T_0^* = G_0^I - T_0^{g*}\right) \\ \hat{Y}_0^* &= \sum_{j \in E} \lambda_j^* Y_j \left(= Y_0 - S_0^* = G_0^O - S_0^{g*}\right)\end{aligned} \quad (5)$$

where $T_0^* = \left(t_{10}^*, ..., t_{m0}^*\right)'$, $T_0^{g*} = \left(t_{10}^{g*}, ..., t_{m0}^{g*}\right)'$, $S_0^* = \left(s_{10}^*, ..., s_{s0}^*\right)'$ and $S_0^{g*} = \left(s_{10}^{g*}, ..., s_{s0}^{g*}\right)'$.

<u>Remark 1.</u> Constraints $\lambda_j d_j = 0$, $j \in E$, are non-linear. Nevertheless, (4) can be solved in practice by reformulating these constraints using Special Ordered Sets (SOS) (Beale and Tomlin, 1970). SOS Type 1 is a set of variables where at most one variable may be nonzero. Therefore, if we remove these constraints from the formulation and define instead a SOS Type 1 for each pair of variables $\{\lambda_j, d_j\}$, $j \in E$, then it is ensured that $\lambda_j$ and $d_j$ cannot be simultaneously positive for DMU$_j$'s, $j \in E$. CPLEX Optimizer (and also LINGO) can solve LP problems with SOS. SOS variables have already been used for solving models like (4) in Ruiz and Sirvent (2016), Aparicio et al. (2017), Cook et al. (2017) and Cook et al. (2018).

<u>Remark 2.</u> Model (4) is non-linear as a consequence of the use of absolute values in the objective function. However, it can be easily linearized by means of a change of variables.

## 3. Benchmarking adjusted to goals oriented either to inputs or to outputs

In practice, goals are often established only for a set of outputs. For instance, in the example examined in this paper, wherein teaching performance of universities is concerned, goals are typically set only for indicators like the rates of progress, drop out and graduation. In these situations, the monitoring of improvements should be made through an evaluation of performance based on output targets that represent best practices which can be implemented with the actual level of resources. In other words, we need to develop an output-oriented version of the benchmarking model formulated in the previous section. A similar reasoning could be made in case of having goals oriented only to inputs. The developments in this section are made in the output oriented case. The extension to input oriented models is straightforward.[6]

Now, we suppose that, at a prior stage to the evaluation, goals $G_0^O = \left(g_{10}^O, \ldots, g_{s0}^O\right)'$ were established for the outputs of a given $DMU_0$, whose actual inputs/outputs in the subsequent period of performance are $(X_0, Y_0)$ (as in the previous section, $(X_j, Y_j)$ also represent the inputs and outputs of all the DMUs in the evaluation period). For the setting of targets in this framework, the key is in defining what is meant by a benchmark. In the output-oriented case, benchmarks should be selected among the points of the PPS having inputs not greater than $X_0$ and outputs that cannot be improved with the actual level of inputs, $X_0$. That is, points $(\hat{X}_0, \hat{Y}_0)$ of the PPS with $\hat{X}_0 \leq X_0$ and whose output bundle $\hat{Y}_0$ belongs to $\partial(P(X_0))$, which is the efficient frontier of $P(X_0) = \{Y / (X_0, Y) \in T\}$. In particular, the output targets $\hat{Y}_0^*$ satisfying the objectives of closeness to both actual performances and goals can be found as the optimal solution of the following model

$$\text{Min} \quad \alpha \left\| Y_0 - \hat{Y}_0 \right\|_1^{\omega_o} + (1-\alpha) \left\| G_0^O - \hat{Y}_0 \right\|_1^{\omega_o}$$
$$\text{s.t.:} \quad \hat{Y}_0 \in \partial(P(X_0)) \tag{6}$$

where $\alpha$, $0 \leq \alpha \leq 1$, is used for the same purposes as in (3) and $\left\| Y_0 - \hat{Y}_0 \right\|_1^{\omega_o} = \sum_{r=1}^{s} \frac{|y_{r0} - \hat{y}_{r0}|}{y_{r0}}$ and $\left\| G_0^O - \hat{Y}_0 \right\|_1^{\omega_o} = \sum_{r=1}^{s} \frac{|g_{r0} - \hat{y}_{r0}|}{y_{r0}}$.

To find an operational formulation of (6) we need to characterize the set $\partial(P(X_0))$.

---

[6] See Aparicio et al. (2017) for an approach based on bilevel linear programming for determining the least distance to the efficient frontier and setting the closest targets that uses an oriented version of the Russell measure.

**Proposition 2.**

$\hat{Y}_0 \in \partial(P(X_0)) \Leftrightarrow X_0 = \sum_{j \in E} \lambda_j X_j + T_0, \ \hat{Y}_0 = \sum_{j \in E} \lambda_j Y_j, \ \sum_{j \in E} \lambda_j = 1,$ for some $T_0 \geq 0_m$, and $\lambda_j \geq 0, \ j \in E,$ satisfying $\lambda_j d_j = 0, \ j \in E,$ and $v_i t_{i0} = 0, \ i = 1, ..., m,$ where $-v'X_j + u'Y_j + u_0 + d_j = 0, \ j \in E, \ v \geq 0_m,$ $u \geq 1_s, u_0 \in \mathbb{R}.$

<u>Proof.</u> It is well-known that $\hat{Y}_0 \in \partial(P(X_0))$ if, and only if, the optimal value of the following problem equals zero: $\text{Max} \left\{ \sum_{r=1}^{s} s_{r0}^+ \ \middle| \ \sum_{j \in E} \lambda_j x_{ij} + t_{i0} = x_{i0}, i = 1, ..., m, \sum_{j \in E} \lambda_j y_{rj} + s_{r0}^+ = \hat{y}_{r0}, r = 1, ..., s, \sum_{j \in E} \lambda_j = 1, \right.$
$\left. \lambda_j \geq 0, j \in E, t_{i0} \geq 0, i = 1, ... m, s_{r0}^+ \geq 0, r = 1, ..., s \right\}.$ And the optimal value of this problem is zero if, and only if, there exists a feasible solution with $s_{r0}^+ = 0, \ r = 1, ..., s, \ t_{i0} \geq 0, i = 1, ..., m,$ satisfying the optimality conditions. That is, if, and only if, there exist $v \geq 0_m, \ u \geq 1_s, \ d_j, \lambda_j \geq 0, \ j \in E, \ u_0 \in \mathbb{R},$ such that $\sum_{j \in E} \lambda_j x_{ij} + t_{i0} = x_{i0}, \quad i = 1, ..., m, \quad \sum_{j \in E} \lambda_j y_{rj} = \hat{y}_{r0}, \quad r = 1, ..., s, \quad \sum_{j \in E} \lambda_j = 1,$
$-v'X_j + u'Y_j + u_0 + d_j = 0, \ j \in E,$ satisfying $\lambda_j d_j = 0, j \in E,$ and $v_i t_{i0} = 0, i = 1, ..., m.$ ∎

Taking into account the proposition above, we can use the following operational formulation of model (6) in order to find the output targets $\hat{Y}_0^*$

$$\text{Min} \quad \alpha \sum_{r=1}^{s} \frac{|s_{r0}|}{y_{r0}} + (1-\alpha) \sum_{r=1}^{s} \frac{|s_{r0}^g|}{y_{r0}}$$

s.t.:

$$\sum_{j \in E} \lambda_j x_{ij} = x_{i0} - t_{i0} \qquad i = 1,...,m$$

$$\sum_{j \in E} \lambda_j y_{rj} = y_{r0} + s_{r0} \qquad r = 1,...,s$$

$$\sum_{j \in E} \lambda_j y_{rj} = g_{r0}^O + s_{r0}^g \qquad r = 1,...,s$$

$$\sum_{j \in E} \lambda_j = 1$$

$$-\sum_{i=1}^{m} v_i x_{ij} + \sum_{r=1}^{s} u_r y_{rj} + u_0 + d_j = 0 \qquad j \in E$$

$$u_r \geq 1 \qquad r = 1,...,s$$
$$\lambda_j d_j = 0 \qquad j \in E$$
$$v_i t_{i0} = 0 \qquad i = 1,...,m$$
$$\lambda_j, d_j \geq 0 \qquad j \in E$$
$$v_i, t_{i0} \geq 0 \qquad i = 1,...,m \qquad (7)$$
$$u_0, s_{r0}, s_{r0}^g \text{ free} \qquad r = 1,...,s$$

<u>Remark 3.</u> As with model (4), the non-linear constraints in (7) can be implemented by using SOS 1 variables when solving that model in practice. Likewise, the absolute values can be avoided by using a change of variables.

For every $\alpha$, the output targets $\hat{Y}_0^*$ can be obtained by using the optimal solution of (7) as follows

$$\hat{Y}_0^* = \sum_{j \in E} \lambda_j^* Y_j \left(= Y_0 - S_0^* = G_0^O - S_0^{g*}\right) \qquad (8)$$

where $S_0^* = \left(s_{10}^*,...,s_{s0}^*\right)'$ and $S_0^{g*} = \left(s_{10}^{g*},...,s_{s0}^{g*}\right)'$.

## 4. Dealing with non-controllable variables

The models developed so far assume that all the inputs and outputs considered for the evaluation of performance are controllable. However, in practice, the levels of some of these variables are exogenously fixed frequently. That is the case, for example, of size and scientific specialization in the evaluation of university performance or the mean income in the demographic area of the unit under evaluation. In this section, we show how the DEA benchmarking adjusted to goals can be readily adapted for the case where at least one of the inputs and/or outputs are non-controllable. In order to do so, we follow the ideas in the approach proposed in Banker and Morey (1986) for dealing with non-discretionary variables.

For the developments, we consider the index sets $I_D$ and $I_{ND}$, which correspond to the controllable and non-controllable inputs, respectively, so that $I = \{1,...,m\} = I_D \cup I_{ND}$. Let $m_D$ and $m_{ND}$ be the cardinalities of both sets, respectively (therefore, $m = m_D + m_{ND}$). Similarly, we have the index sets for the outputs $O_D$ and $O_{ND}$, $O = \{1,...,s\} = O_D \cup O_{ND}$, so that $s = s_D + s_{ND}$, $s_D$ and $s_{ND}$ being the corresponding cardinalities. Let $X_0^D$ and $X_0^{ND}$ be the vectors of controllable and non-controllable inputs of a given DMU$_0$, $Y_0^D$ and $Y_0^{ND}$ being the ones corresponding to controllable and non-controllable outputs. For the benchmarking of DMU$_0$, we can consider the set $T(X_0^{ND}, Y_0^{ND}) = \{(X^D, Y^D) / (X^D, X_0^{ND}, Y^D, Y_0^{ND}) \in T\}$, where T is the production possibility set associated with the whole set of variables. $T(X_0^{ND}, Y_0^{ND})$ is actually the vector set corresponding to controllable inputs and outputs which are compatible with $(X_0^{ND}, Y_0^{ND})$. In a similar manner as in the previous section, the benchmarks for a given DMU$_0$ with actual inputs/outputs $(X_0^D, X_0^{ND}, Y_0^D, Y_0^{ND})$ should be selected among the points of the PPS with controllable inputs/outputs $(\hat{X}_0^D, \hat{Y}_0^D)$ that cannot be improved with the actual levels of non-controllable inputs and outputs $(X_0^{ND}, Y_0^{ND})$. That is, points $(\hat{X}_0^D, \hat{X}_0^{ND}, \hat{Y}_0^D, \hat{Y}_0^{ND})$ of the PPS with $\hat{X}_0^{ND} \leq X_0^{ND}$ and $\hat{Y}_0^{ND} \geq Y_0^{ND}$ whose controllable input-output bundle $(\hat{X}_0^D, \hat{Y}_0^D)$ belongs to $\partial(T(X_0^{ND}, Y_0^{ND}))$, the efficient frontier of $T(X_0^{ND}, Y_0^{ND})$. If we denote by $G_0^{I_D}$ and $G_0^{O_D}$ the vectors consisting of the goals established for the controllable inputs and outputs, respectively, the model below provides targets $(\hat{X}_0^{D*}, \hat{Y}_0^{D*})$ according to the approach followed in this paper:

$$\begin{aligned}
&\text{Min} \quad \alpha \left\|(X_0^D, Y_0^D) - (\hat{X}_0^D, \hat{Y}_0^D)\right\|_1^\omega + (1-\alpha) \left\|(G_0^{I_D}, G_0^{O_D}) - (\hat{X}_0^D, \hat{Y}_0^D)\right\|_1^\omega \\
&\text{s.t.:} \\
&\qquad (\hat{X}_0^D, \hat{Y}_0^D) \in \partial\left(T(X_0^{ND}, Y_0^{ND})\right)
\end{aligned} \qquad (9)$$

for a given $\alpha$, $0 \leq \alpha \leq 1$, where $\left\|\left(X_0^D, Y_0^D\right)-\left(\hat{X}_0^D, \hat{Y}_0^D\right)\right\|_1^\omega = \sum_{i \in I_D} \frac{|x_{i0} - \hat{x}_{i0}|}{x_{i0}} + \sum_{r \in O_D} \frac{|y_{r0} - \hat{y}_{r0}|}{y_{r0}}$ and

$\left\|\left(G_0^{ID}, G_0^{OD}\right)-\left(\hat{X}_0^D, \hat{Y}_0^D\right)\right\|_1^\omega = \sum_{i \in I_D} \frac{|g_{i0} - \hat{x}_{i0}|}{x_{i0}} + \sum_{r \in O_D} \frac{|g_{r0} - \hat{y}_{r0}|}{y_{r0}}$.

The following results will allow us to derive an operational formulation of model (9):

**Lemma.**

$\left(\hat{X}_0^D, \hat{Y}_0^D\right) \in \partial\left(T\left(X_0^{ND}, Y_0^{ND}\right)\right)$ if, and only if, the optimal value of the following model is zero:

$\text{Max}\left\{\sum_{i \in I_D} t_{i0} + \sum_{r \in O_D} s_{r0} \,\middle|\, \sum_{j \in E} \lambda_j x_{ij} + t_{i0} = \hat{x}_{i0}, i \in I_D, \sum_{j \in E} \lambda_j x_{ij} \leq x_{i0}, i \in I_{ND}, \sum_{j \in E} \lambda_j y_{rj} - s_{r0} = \hat{y}_{r0}, r \in O_D,\right.$

$\left.\sum_{j \in E} \lambda_j y_{rj} \geq y_{r0}, r \in O_{ND}, \sum_{j \in E} \lambda_j = 1, t_{i0} \geq 0, i \in I_D, s_{r0} \geq 0, r \in O_D, \lambda_j \geq 0, j \in E\right\}$.

The following proposition provides a characterization of $\partial\left(T\left(X_0^{ND}, Y_0^{ND}\right)\right)$

**Proposition 3.**

$\left(\hat{X}_0^D, \hat{Y}_0^D\right) \in \partial\left(P(X_0^{ND}, Y_0^{ND})\right) \Leftrightarrow \hat{X}_0^D = \sum_{j \in E} \lambda_j X_j^D, \quad X_0^{ND} = \sum_{j \in E} \lambda_j X_j^{ND} + T_0^{ND}, \quad \hat{Y}_0^D = \sum_{j \in E} \lambda_j Y_j^D,$

$Y_0^{ND} = \sum_{j \in E} \lambda_j Y_j^{ND} - S_0^{ND}, \sum_{j \in E} \lambda_j = 1,$ for some $T_0^{ND} \geq 0_{m_{ND}}$, $S_0^{ND} \geq 0_{s_{ND}}$ and $\lambda_j \geq 0$, $j \in E$, satisfying

$\lambda_j d_j = 0, \; j \in E, \; v_i t_{i0} = 0, \; i \in I_{ND},$ and $u_r s_{r0} = 0, \; r \in O_{ND}$, where $-v'X_j + u'Y_j + u_0 + d_j = 0, \; j \in E,$

$v = \left(v^D, v^{ND}\right), \; v^D \geq 1_{m_D}, \; v^{ND} \geq 0_{m_{ND}}, \; u = \left(u^D, u^{ND}\right), \; u^D \geq 1_{s_D}, \; u^{ND} \geq 0_{s_{ND}}, \; u_0 \in \mathbb{R}.$

Proof. The proof of this proposition is similar to that of Proposition 2. That is, the necessary and sufficient conditions for $\left(\hat{X}_0^D, \hat{Y}_0^D\right)$ to belong to $\partial\left(T\left(X_0^{ND}, Y_0^{ND}\right)\right)$ result from the optimality conditions associated with the linear problem in the Lemma when the optimal value equals zero.

Proposition 3 leads us to the following reformulation of model (9):

$$\text{Min} \quad \alpha\left(\sum_{i\in I_D}\frac{|t_{i0}|}{x_{i0}}+\sum_{r\in O_D}\frac{|s_{r0}|}{y_{r0}}\right)+(1-\alpha)\left(\sum_{i\in I_D}\frac{|t_{i0}^g|}{x_{i0}}+\sum_{r\in O_D}\frac{|s_{r0}^g|}{y_{r0}}\right)$$

s.t.:

$$\sum_{j\in E}\lambda_j x_{ij}=x_{i0}-t_{i0} \qquad i=1,\ldots,m$$

$$\sum_{j\in E}\lambda_j y_{rj}=y_{r0}+s_{r0} \qquad r=1,\ldots,s$$

$$\sum_{j\in E}\lambda_j x_{ij}=g_{i0}^I-t_{i0}^g \qquad i\in I_D$$

$$\sum_{j\in E}\lambda_j y_{rj}=g_{r0}^O+s_{r0}^g \qquad r\in O_D$$

$$\sum_{j\in E}\lambda_j=1$$

$$-\sum_{i=1}^m v_i x_{ij}+\sum_{r=1}^s u_r y_{rj}+u_0+d_j=0 \qquad j\in E$$

$$\lambda_j d_j=0 \qquad j\in E$$

$$v_i t_{i0}=0 \qquad i\in I_{ND}$$

$$u_r s_{r0}=0 \qquad r\in O_{ND}$$

$$v_i\geq 1 \qquad i\in I_D$$

$$u_r\geq 1 \qquad r\in O_D$$

$$v_i\geq 0, t_{i0}\geq 0 \qquad i\in I_{ND} \qquad (10)$$

$$u_r\geq 0, s_{r0}\geq 0 \qquad r\in O_{ND}$$

$$d_j,\lambda_j\geq 0 \qquad j\in E$$

$$t_{i0}, t_{i0}^g \text{ free} \qquad i\in I_D$$

$$u_0, s_{r0}, s_{r0}^g \text{ free} \qquad r\in O_D$$

<u>Remark 4.</u> If $I_{ND}$ or $O_{ND}$ are $\varnothing$, that is, if there is either no non-controllable input or no non-controllable output, then we just need to remove from (10) the constraints and the variables associated with the corresponding index set.

<u>Remark 5.</u> In case of having goals only for outputs, an oriented model for performance evaluation can be obtained by simply setting $I_D=\varnothing$ in (10) (note that output-oriented models implicitly treats inputs as non-controllable; see Thanassoulis et al. (2008) for discussions). In case of having only input oriented goals, we should set $O_D=\varnothing$.

## 5. Illustrative example

For purposes of illustration only, in this section we apply the proposed approach to the evaluation of teaching performance of public Spanish universities, within the context of the implementation of plans for improvements. As said in the introduction, public Spanish universities often implement plans for improvements by setting goal levels for a number of index indicators, in particular regarding teaching performance. Regional governments and university managers set the goals and an evaluation of performance is made in the next academic year in order to assess the degree of achievements of goals. In fact, these can be seen as incentive plans, because part of the financing of public universities (up to 10% of base financing) is linked to performance relative to goals. In the analysis that is carried out here, we suppose that stakeholders, in the academic year 2013-14, established goals for each university regarding the results of their teaching activities, and that performance is to be evaluated relative to the targets determined by the benchmarking model with the actual data corresponding to 2014-15.

The public Spanish universities may be seen as a set of homogeneous DMUs that undertake similar activities and produce comparable results regarding, in particular, teaching performance, so that a common set of outputs can be defined for their evaluation. Specifically, for the selection of the outputs, we take into consideration variables that represent aspects of performance explicitly mentioned as requirements in the section "Expected Results" of the Royal Order RD 1393/2007, such as graduation, drop out and progress of students. As for the inputs, we select two variables that account for human and physical capital: academic staff and expenditures. In addition, the size of the universities (proxied by the number of students) is incorporated into the analysis as a non-controllable factor. As recognized in the related literature, size has a significant impact on the performance of European universities (see, Daraio et al. (2015), which includes an exhaustive literature review). The variables considered are defined as explained below. The data correspond to the academic year 2014-15, which is adopted as the reference year (see, for example, Agasisti and Dal Bianco (2009), which follows this approach for performance evaluation of universities, by selecting variables of this type).

OUTPUTS
- GRADUATION (GRAD): Total number of students that complete the programme of studies within the planned time.
- RETENTION (RET): Total number of students enrolled for the first time in the academic year 2012-13 that keep enrolled at the university.

- PROGRESS (PROG): Total number of passed credits[7].

INPUTS
- ACADEMIC STAFF (ASTF): Full-time equivalent academic staff.
- EXPENDITURE (EXP): This input exactly accounts for expenditure on goods and services after the item corresponding to works carried out by other (external) companies has been removed. EXP thus reflects the budgetary effort made by the universities in the delivery of their activities.

SIZE, which is computed as the total number of students enrolled in the academic year 2014-15, is treated as a non-discretionary input.

Data for these variables have been taken from the corresponding report by the Conference of Rectors of the Spanish Universities (CRUE). The sample consists of 38 (out of 48) public Spanish universities. Table 1 reports a descriptive summary.

Table 1. descriptive summary

For the benchmarking and the setting of targets, we use an output-oriented model that includes as inputs and outputs the ones listed above and incorporates the size of universities as a non-controllable factor. An initial DEA analysis reveals 18 universities as technically efficient. Table 2 reports the results derived from the benchmarking model regarding some universities that have been chosen as representative cases: the universities of Cádiz (UCA), Alicante (UA), Huelva (UHU), Autónoma de Barcelona (UAB), Zaragoza (UZA) and Santiago de Compostela (USC). For each of these universities, the table records its actual data in the academic year 2014-15, the targets provided by the model for a grid of α values between 0 and 1 and its goals. As said before, we suppose that, in the academic year 2013-14, governments and university managers set output goals for each university as specified in the rows "Goals" of Table 2. We point out that these goals are quite similar across universities (if we take into account their sizes), which means that some overall management objectives are pursued regarding each of the outputs. The small differences among them only seek to consider to some extent individual circumstances.

Firstly, it should be pointed out the fact that the goals that were set in the academic year 2013-14 are attainable only in the cases of UZA and USC, while the points consisting of actual inputs, size (in 2014-15) and goals are outside the PPS in the cases of UCA, UA, UHU and UAB. Thus, the latter four universities could argue that the goals established for them are unachievable taking into account their actual levels of ASTAFF and EXPEND as well as their SIZE in the academic year 2014-15, whereas in the case UZA and USC these could be deemed as little too unambitious goals.

---

[7] Credit is the unit of measurement of the academic load of the subject of a programme.

For the evaluation of UCA, the closest targets (those associated with high values of α) point to a weakness in RET. However, when the benchmarking is adjusted to the goals (see the rows with lower values of α), while setting attainable targets, we can see that there is also room for improvement in GRAD. Something similar occurs to UA. In this case, the adjustment to the goals entails setting a closer target for RET than those associated with higher α's, in exchange for setting targets somewhat more demanding in GRAD and PROG. Eventually, the evaluation relative to these targets reveals that UA needs to improve in all of the three outputs. UHU was rated as efficient. However, when the target setting is adjusted to the goals, it is shown that there is room for improvement, particularly in GRAD. In contrast, UAB, which was also efficient, shows a good performance in all of the outputs in all of the scenarios of benchmarking, in spite of having set for this university the most demanding goals.

Actual outputs of UZA show that this university is performing at the level of the goals. However, as said before, these are inside the PPS and can be deemed as little too unambitious goals. In fact, the model sets targets more demanding than the goals particularly for RET (even for low values of α), thus suggesting that UZA should improve its performance regarding the abandonment of students. The situation with USC is similar. There is room for improvement for USC in all of the three outputs, in spite of performing not far from the goals. This is particularly noticeable when it is evaluated against the targets associated with lower α's, which are more demanding (they are at the goal level in PROG and are quite higher than the goals in GRAD and RET).

## 6. Conclusions

We have shown how DEA can be extended to the evaluation of performance within the framework of the implementation of plans for improvements that set management goals. DEA benchmarking models have been formulated which allow us to set targets adjusted to the goals, in addition to being attainable and representing best practices. The approach proposed is developed under the basic philosophy of DEA as ex post tool for evaluating performances within a process of monitoring and control. However, as has been pointed out in the related literature, the benchmarking for decision making units goes beyond a purely monitoring process, and includes a component of future planning as a natural step (see Thanassoulis et al. (2008), Yang et al. (2009) and Stewart (2010)). These authors have also shown how DEA can be adapted for use with those purposes. In this sense, it should be highlighted that the models proposed in this paper can be utilized for setting future targets as well, thus becoming an ex ante planning tool. Specifically, the same formulations developed here can be solved with data of actual performances in a given period of time (when goals are established) in order to set targets for a subsequent period of performance.

As a future direction, the extension of the approach proposed to deal with the situation in which the DMUs are organized into groups of units whose member experience similar circumstances (for example, universities of a given region that operate under the same regulatory framework) becomes an interesting line of research. The idea would be to develop a within-group common benchmarking framework in line with that proposed in Cook et al. (2017), while at the same time incorporating the information on management goals. We would like also to investigate more deeply the case where some of the inputs and/or outputs are non-controllable. The model proposed in section 4 has been developed following the ideas in Banker and Morey (1986). However, some authors have pointed out some weaknesses in the standard approach proposed by these authors to deal with non-discretionary inputs and outputs. For this reason, new formulations that extend the basic model (10) should be developed that allow us address these issues more properly within the context of the benchmarking adjusted to goals.

**Acknowledgments**

This research has been supported through Grant MTM2016-76530-R (AEI/FEDER, UE).

|  | SIZE | ASTAFF | EXPEND | GRAD | RET | PROG |
|---|---|---|---|---|---|---|
| Mean | 20362.05 | 1639.87 | 15513793.90 | 1845.18 | 3801.42 | 795103.24 |
| Standard Dev. | 12105.02 | 949.43 | 9828600.31 | 1057.09 | 2197.37 | 466860.90 |
| Minimum | 3735 | 371.88 | 2514366.67 | 399 | 755 | 140645.5 |
| Maximum | 55662 | 3855.13 | 46471378.25 | 4804 | 9442 | 1983413 |

Table 1. Descriptive summary

Table 2. Targets for different values of α

| UCA | | SIZE | ASTAFF | EXPEND | GRAD | RET | PROG |
|---|---|---|---|---|---|---|---|
| | Actual | 18605 | 1548.5 | 9908074.26 | 1169 | 3117 | 738346 |
| | α=1 | | 889.92 | 9908074.26 | 1169 | 3797.79 | 738329.69 |
| | α=0.8 | | 890.15 | 9908074.26 | 1169.07 | 3797.74 | 738346 |
| | α=0.6 | | 890.15 | 9908074.26 | 1169.07 | 3797.74 | 738346 |
| | α=0.5 | | 1051.31 | 9908074.26 | 1593.77 | 3697.96 | 747379.10 |
| | α=0.4 | | 1219.42 | 9908074.26 | 1700 | 3644.14 | 742409.82 |
| | α=0.2 | | 1219.42 | 9908074.26 | 1700 | 3644.14 | 742409.82 |
| | α=0 | | 1219.42 | 9908074.26 | 1700 | 3644.14 | 742409.82 |
| | Goals* | | | | **1700** | **3700** | **765000** |
| UA | | SIZE | ASTAFF | EXPEND | GRAD | RET | PROG |
| | Actual | 24322 | 2036.25 | 15204774.99 | 2378 | 4234 | 948324 |
| | α=1 | | 1767.55 | 15204774.99 | 2378 | 4885.41 | 950313.42 |
| | α=0.8 | | 1767.55 | 15204774.99 | 2378 | 4885.41 | 950313.42 |
| | α=0.6 | | 1767.55 | 15204774.99 | 2378 | 4885.41 | 950313.42 |
| | α=0.5 | | 1767.55 | 15204774.99 | 2378 | 4885.41 | 950313.42 |
| | α=0.4 | | 1897.66 | 15204774.99 | 2557.51 | 4690.45 | 1001429.16 |
| | α=0.2 | | 1897.66 | 15204774.99 | 2557.51 | 4690.45 | 1001429.16 |
| | α=0 | | 1897.66 | 15204774.99 | 2557.51 | 4690.45 | 1001429.16 |
| | Goals* | | | | **2800** | **4500** | **940000** |
| UHU† | | SIZE | ASTAFF | EXPEND | GRAD | RET | PROG |
| | Actual | 10601 | 706.63 | 5903210.21 | 525 | 2306 | 378848.5 |
| | α=1 | | 706.63 | 5903210.21 | 525 | 2306 | 378848.5 |
| | α=0.8 | | 706.63 | 5903210.21 | 525 | 2306 | 378848.5 |
| | α=0.6 | | 706.63 | 5903210.21 | 525 | 2306 | 378848.5 |
| | α=0.5 | | 706.63 | 5903210.21 | 525 | 2306 | 378848.5 |
| | α=0.4 | | 706.63 | 5903210.21 | 848.51 | 2162.60 | 410000 |
| | α=0.2 | | 706.63 | 5903210.21 | 850 | 2161.94 | 410143.33 |
| | α=0 | | 706.63 | 5903210.21 | 850 | 2161.94 | 410143.33 |
| | Goals* | | | | **850** | **2500** | **410000** |
| UAB† | | SIZE | ASTAFF | EXPEND | GRAD | RET | PROG |
| | Actual | 26983 | 2499.375 | 32279330.43 | 3573 | 5628 | 1238228 |
| | α=1 | | 2499.38 | 32279330.43 | 3573 | 5628 | 1238228 |
| | α=0.8 | | 2499.38 | 32279330.43 | 3573 | 5628 | 1238228 |
| | α=0.6 | | 2499.38 | 32279330.43 | 3573 | 5628 | 1238228 |
| | α=0.5 | | 2499.38 | 32279330.43 | 3573 | 5628 | 1238228 |
| | α=0.4 | | 2460.99 | 29529852.92 | 3300 | 5785.69 | 1240000 |
| | α=0.2 | | 2460.99 | 29529852.92 | 3300 | 5785.69 | 1240000 |
| | α=0 | | 2460.99 | 29529852.92 | 3300 | 5785.69 | 1240000 |
| | Goals* | | | | **3300** | **6000** | **1240000** |

| UZA | | SIZE | ASTAFF | EXPEND | GRAD | RET | PROG |
|---|---|---|---|---|---|---|---|
| | Actual | 27054 | 2834.5 | 22531331 | 2512 | 5179 | 1115003 |
| | α=1 | | 2216.87479 | 22531331 | 2621.42 | 5746.42 | 1115003 |
| | α=0.8 | | 2216.87479 | 22531331 | 2621.42 | 5746.42 | 1115003 |
| | α=0.6 | | 2216.87479 | 22531331 | 2621.42 | 5746.42 | 1115003 |
| | α=0.5 | | 2216.87479 | 22531331 | 2621.42 | 5746.42 | 1115003 |
| | α=0.4 | | 2216.87479 | 22531331 | 2621.42 | 5746.42 | 1115003 |
| | α=0.2 | | 2233.57359 | 22531331 | 2640.80 | 5739.97 | 1130000 |
| | α=0 | | 2233.57359 | 22531331 | 2640.80 | 5739.97 | 1130000 |
| | Goals | | | | **2600** | **5300** | **1130000** |
| USC | | SIZE | ASTAFF | EXPEND | GRAD | RET | PROG |
| | Actual | 20876 | 1889.375 | 17513725.15 | 1731 | 3211 | 820284 |
| | α=1 | | 1724.82 | 17513725.15 | 2014.59 | 4574.09 | 842587.07 |
| | α=0.8 | | 1724.82 | 17513725.15 | 2014.59 | 4574.09 | 842587.07 |
| | α=0.6 | | 1724.82 | 17513725.15 | 2014.59 | 4574.09 | 842587.07 |
| | α=0.5 | | 1724.82 | 17513725.15 | 2014.59 | 4574.09 | 842587.07 |
| | α=0.4 | | 1724.82 | 17513725.15 | 2014.59 | 4574.09 | 842587.07 |
| | α=0.2 | | 1744.21 | 17513725.15 | 2037.09 | 4566.60 | 860000 |
| | α=0 | | 1744.21 | 17513725.15 | 2037.09 | 4566.60 | 860000 |
| | Goals | | | | **1800** | **3500** | **860000** |

* Goals unachievable for the actual level of resources and size

† Universities DEA efficient